\documentclass{tran-l}

\newtheorem{theorem}{Theorem}[subsection]

\newtheorem{lemma}[theorem]{Lemma}

\theoremstyle{definition}
\newtheorem{definition}[theorem]{Definition}
\theoremstyle{remark}

\numberwithin{equation}{subsection}

\begin{document}

\title{Closed-Form Trajectories of the Unicycle.}

\author{Andr\'{e}s Garc\'{i}a ,Osvaldo Agamennoni and Jos\'{e} Figueroa}

\address{Department of Electrical Engineering and Computers, Universidad Nacional del Sur, Bah\'{i}a Blanca, Buenos Aires,
 Argentina}

\email{agarcia@uns.edu.ar}


\begin{abstract}
This paper presents a methodology to stabilize some kind of Nonlinear Control system known as Driftless, utilizing the concept of \textit{Pseudo-Kinetic Energy} introduced in this work. Once this controller is applied to the Unicycle-type robot, stability is guaranteed with the salient property that the structure of the controller allows to solve in closed-form the trajectories of the vehicle.

While the proposed controller only ensures stability (not asymptotic stability) the obtained closed-form solutions will show a path to obtain in closed form the solutions for the general control problem of the unicycle.

Some conclusions and future directions for research are also depicted.

\end{abstract}


\maketitle

\section{Introduction}

The stabilization of Driftless affine systems is a very important topic both in practical cases and  theoretical analysis (Non holonomic vehicles, robotic systems in chain
form,etc.), because of that, they have received a lot of attention in past years (see for instance \cite{Morin03} and \cite{Lin96} with the references trough them),

While all the available techniques for stabilization of Nonlinear Driftless systems consider time varying open-loop control laws (\cite{Duleba98}, \cite{Murray93}) or smooth controllers using Lyapunov Functions (see \cite{Aicardi94}); up to the authors knowledge the implementation of a controller for Driftless systems with the possibility of Closed-Form solutions of stable trajectories was never done in the past.

In this way, the goal in this paper will be to derive a general continuous smooth stabilizing control law for Driftless systems, introducing the concept of
\textit{Pseudo-Kinetic Energy}. This control law will be applied to the well known Unicycle-like robot, yielding Closed-Form solutions for the trajectories of the wheel baseline center coordinates and the attitude angle via the solution of a second order Bessel's equation.

Once the closed-form solutions are obtained, the impossibility of an asymptotic behavior for one of the wheel baseline center coordinates will be established in accordance with the Brockett's theorem (see \cite{Brockett83}). However, as it will be shown, if the same controller is used in such a way that the attitude is unstable (the robot continues rotating), then the wheel baseline center will be in fact asymptotically stable.

This conclusion will guide an extension for this controller via the addition of a dissipative term which produces some sort of "shaking" behavior of the attitude, endowing the entire Kinematic Model of the vehicle with asymptotic stability. It turns out that the extended controller is time-varying preventing collision with Brockett's condition.

It is worth mentioning that controllability is not a condition to be checked in the application of this methodology, even though the controllable nature of the system is not discussed in this work, and the contrast with cases where the Lie Algebra based test fails while the system remains controllable is not focused in this paper.

This paper is organized as follows: Section \ref{El problema considerado} introduces formally the class of systems considered in the paper, Section \ref{Estabilizacion de Driftless} presents general theorems for stabilization of Driftless systems, Section \ref{Uniciclo} apply the previous theory to the Unicycle and finally Section \ref{Uniciclo asintotico} extends the analysis for asymptotic stability of the Unicycle.

\section{The Problem Considered and Motivating Ideas}\label{El problema considerado}

In this section, the problem considered in this paper and some preliminary ideas related with the stability criterion used to design stabilizing smooth controllers will be established.

The kind of Kinematic Models considered in this paper are the well known "Driftless" systems:

\begin{equation}\label{Definicion de los sistemas a usar}
\dot{q}(t)=S(q)\cdot \tilde{u}(q)\text{,}
\end{equation}

where $q(t) \in \Re^{n\times 1}$ are the state space variables, $\tilde u(q) \in \Re ^{k \times 1}$  is the admissible control set which is at least $C^1(M)$, where $M$ could be any manifold configuration for the state space $q(t)$.

The challenge is to develop a systematic technique with two main characteristics:

\begin{itemize}
\item{Stabilization of the system (\ref{Definicion de los sistemas a usar}) in Lyapunov sense.}
\item{The possibility of addition of dissipative terms in order to make the system asymptotically stable.}
\end{itemize}

The first item will be reached introducing the concept of Pseudo-Kinetic Energy and showing that such a definition allows obtaining a Lyapunov function which ensures stability for Driftless systems.  Realize at this point that asymptotic stability is not yet addressed.

A special Lyapunov function is obtained by introducing the concept of Pseudo-Kinetic Energy for a subclass of Driftless systems.  It is showed that the Unicycle is in this category. Moreover, the technique applied to the Unicycle yields a stable system and also allows to have closed-form solutions for the trajectories of such a system.

Analyzing these closed-form solutions is concluded the impossibility of asymptotic stability with this controller.  The problem could be solved with the inclusion of a dissipative term. In fact, by defining a general dissipative term for the attitude angle of the Unicycle, it will be shown in an analytic way that the Kinematic Model will turn to be asymptotically stable.

\section{A Stabilizing Controller for Driftless Systems}\label{Estabilizacion de Driftless}

From now on, stability will refers to stability in the sense of Lyapunov.  So, a stabilizing controller should satisfy the Lyapunov theorem which yields a sufficient condition for a Nonlinear system to be stable (see \cite{Sastry99}, pp. 188-189). Realize that the concept of l2-norm plays a fundamental role because it is a scalar quantity that allows to measure how near is the system to a stable equilibrium point.

In this way, since the focus is to drive the system to an equilibrium point ($\dot q = 0$), and recalling that the Barbalat's lemma predicts that if an integrand is going to zero anytime the integral is bounded, and with afore mentioned the idea of the l2-norm in mind; it is reasonable to consider the following integral\footnote{Also consider the ideas in the paper \cite{Byrnes95}}:

\begin{equation}\label{Ideas Primarias}
\int_{0}^{t} \langle \dot{q}(\xi ),\dot{q}(\xi )\rangle \cdot d\xi=\text{Bounded},
\end{equation}

where $\left\langle .,.\right\rangle $ means the inner product and $\left\langle \dot{q}(\xi ),\dot{q}(\xi )\right\rangle =\left\Vert \dot{q}(t)\right\Vert ^{2}$.

The machinery in what follows will derive a methodology to design a Control law $\widetilde{u}(q)$ for the Driftless Kinematic Model introduced in equation (\ref{Definicion de los sistemas a usar}). In this way, the first result is an extension of the well known Barbalat's lemma - all the proofs are given in the appendix:

\begin{lemma}\label{Barbalat's lemma extension}
Let $\varphi :\Re
\rightarrow \Re $ be a positive definite function. Suppose that
$\lim_{t\rightarrow \infty}\int\limits_{0}^{t}\varphi (\tau )\cdot
d\tau $ exits and is finite. Then, $\varphi (t)\rightarrow 0$ as
$t\rightarrow \infty .$
\end{lemma}

This lemma shows that one way to guarantee $\dot{q}%
(t)\rightarrow 0$ when $t\rightarrow +\infty $ is to ensure that the integral in (\ref{Ideas Primarias}) when $t\rightarrow +\infty $ is bounded and this motivates the following result:

\begin{lemma}[Stability Criterion]\label{Nuestra definicion de estabilidad}
Given a Control System of the form $\dot{q}(t)=S(q(t))\cdot \tilde{u}(q)$, where $S(.,.)$ is a smooth collection of vector fields which belong to $\Re ^{n\times k}$ with $q=0$ an equilibrium point. Then, this system is stable if and only if the quantity
$E_{c}^{*}=\int_{0}^{\infty }\left\langle \dot{q}(\xi
),\dot{q}(\xi )\right\rangle \cdot d\xi $ called Pseudo-Kinetic
Energy is bounded.
\end{lemma}

Next section is using this Lemma for a subclass of the more general Driftless one, providing a systematic methodology to obtain stabilizing controllers for such a systems.

\section{Stabilizing the Unicycle}\label{Uniciclo}

Based on the Lemma \ref{Nuestra definicion de estabilidad} in previous section, the main result in this paper for a special subclass of driftless systems is presented:

\begin{theorem}[Controller Calculations]\label{Controlador}
Given a Driftless Control system to the form (\ref{Definicion de los sistemas a usar}), such that $\|S_{i}\|=1,\quad i=1,2,\ldots,k$ and $S_{i} \cdot S_{j}=0,\quad i\neq j$ where $\| . \|$ is the $l_{2}$-norm - see \cite{Horn85} for details on matrix norms- and $S=[S_{1},S_{2},\ldots,S_{k}]$.

Then the controller $u_{i}(q)=\rho \cdot q(t)' \cdot S_{i},\quad \rho<0$, stabilize the system $\dot{q}(t)=S(q(t))\cdot \tilde{u}(q)$.
\end{theorem}

Applying this result to the Unicycle, it will be showed that this kind of controller allows to find \textit{Closed-Form solutions for the trajectories} of the system where a posterior modification via "dissipative" terms will produce asymptotic stability.

\subsection{Closed Form Trajectories}\label{Trayectorias en forma cerrada}

Considering the Kinematic Model of the Unicycle-see
\cite{Sastry99}\ pp.529:

\begin{equation}\label{Unicycle}
\dot{q}(t)=\underset{S(q)}{\underbrace{\left[
\begin{array}{cc}
\cos (\theta ) & 0 \\
\sin (\theta ) & 0 \\
0 & 1%
\end{array}
\right] }}\cdot \tilde{u}(q)\text{,}
\end{equation}

where $q(t)=\left[
\begin{array}{ccc}
x_{c}(t) & y_{c}(t) & \theta (t)%
\end{array}%
\right]'$ with $X=\left[
\begin{array}{cc}
x_{c}(t) & y_{c}(t)
\end{array}
\right]'$ the position of the wheel baseline center of the robot and $\theta
(t)$ the attitude of it.

Clearly this system has the property:

\begin{equation*}
\begin{cases}
S_{1}=[\cos(\theta) \quad \sin(\theta) \quad 0]'\\
S_{2}=[0 \quad 0 \quad 1]'\\
\|S_{1}\|=1\\
\|S_{2}\|=1
\end{cases}
\end{equation*}

Then applying the controller in Theorem \ref{Controlador}, we
have:

\begin{equation*}
\tilde{u}(q)=\left[
\begin{array}{c}
\rho \cdot X^{\prime }\cdot \left[
\begin{array}{c}
\cos (\theta ) \\
\sin (\theta ) \\
\end{array}
\right] \\
\rho \cdot \theta
\end{array}
\right] \text{,}
\end{equation*}

where $\rho<0$:

\begin{equation}\label{Separacion de las dinamicas}
\begin{cases}
\left[
\begin{array}{c}
\dot{x}_{c}(t) \\
\dot{y}_{c}(t)
\end{array}
\right] =\rho \cdot \left[
\begin{array}{c}
\cos (\theta ) \\
\sin (\theta )
\end{array}
\right] \cdot \left( \cos (\theta )\cdot xc+\sin (\theta )\cdot yc\right) \\
\dot{\theta}(t)=\rho \cdot \theta (t)
\end{cases}.
\end{equation}

Realize that the dynamic of the attitude ($\theta (t)$) can be decoupled from the rest and it is clearly \emph{asymptotically stable (globally) to zero}. In this way, the concerned in what follows will be the position of the wheel baseline center of the robot $(x_{c}(t),\quad y_{c}(t))$:

\begin{equation*}
\underbrace{
\begin{bmatrix}
\dot{x}_{c}(t)\\
\dot{y}_{c}(t)
\end{bmatrix}}_{\dot{X}(t)}=
\rho \cdot
\begin{bmatrix}
\cos (\theta ) \\
\sin (\theta )
\end{bmatrix}
\cdot (\cos(\theta) \cdot xc + \sin(\theta)\cdot yc)
\end{equation*}

This subsystem can be transformed conveniently via a state space transformation as follows:

\begin{equation}\label{Dinamica del centro de masa}
\underbrace{
\begin{bmatrix}
\dot{x}_{c}(t) \\
\dot{y}_{c}(t)
\end{bmatrix}}_{\dot{X}(t)}
=\rho \cdot \underbrace{
\begin{bmatrix}
\cos (\theta)\\
\sin (\theta)
\end{bmatrix}
\cdot
\begin{bmatrix}
\cos (\theta ) & \sin (\theta)
\end{bmatrix}}_{Y(\theta)} \cdot
\begin{bmatrix}
x_{c}(t) \\
y_{c}(t)
\end{bmatrix}
\end{equation}

Moreover, it is convenient to write $Y(\theta )$ as follows:

\begin{equation*}
Y(\theta )=
\underbrace{
\begin{bmatrix}
\cos(\theta) & - \sin(\theta)\\
\sin(\theta) &  \cos (\theta)
\end{bmatrix}}_{R(\theta)}
\cdot \underbrace{
\begin{bmatrix}
1 & 0 \\
0 & 0
\end{bmatrix}}_{\Delta}
\cdot R(\theta )^{\prime}.
\end{equation*}

It is important to notice that $R(\theta (t))$ is a rotation
matrix about $\theta $ and because of that is always invertible
and smooth. In this way let's consider the following change of
coordinates:

\begin{definition}[Change of Coordinates]\label{Cambio de coordenadas}
$X(t)=R(\theta )\cdot Z(t)$.
\end{definition}

Notice that this transformation is a Lyapunov transformation (see \cite{Rugh96}, pp. 107-108 for details) and because of that preserves the stability properties.

Applying this transformation to (\ref{Dinamica del centro de masa}) with the definition above:

\begin{equation}\label{Dinamica de X(t)}
\dot{X}(t)=\rho \cdot R(\theta)' \cdot \Delta \cdot \underbrace{R(\theta )^{\prime}\cdot \dot{X}(t)}_{Z(t)}
\end{equation}

On the other hand, it is possible to derive an expression for $\dot{Z(t)}$:

\begin{equation*}
\dot{Z}(t)=\dot{R}(\theta )^{\prime }\cdot X(t)+R(\theta )^{\prime
}\cdot \dot{X}(t)
\end{equation*}

Replacing $\dot{X(t)}$ from equation (\ref{Dinamica de X(t)}):

\begin{equation*}
\dot{Z}(t)=\dot{R}(\theta )^{\prime }\cdot X(t)+\rho \cdot \Delta \cdot Z(q)
\end{equation*}

Notice that:

\begin{equation*}
\dot{R}(\theta)^{\prime }=-\dot{\theta}(t)\cdot \underset{R(\theta -\frac{\pi }{2})^{\prime }}{\underbrace{\left[
\begin{array}{cc}
-\sin (\theta) & \cos (\theta) \\
-\cos (\theta) & -\sin (\theta)
\end{array}
\right]}}
\end{equation*}

then:

\begin{gather*}
\dot{Z}(t)=-\dot{\theta}(t)\cdot R(\theta -\frac{\pi }{2})^{\prime
}\cdot R(\theta )\cdot Z(t)+\rho \cdot \Delta \cdot Z(q)\Leftrightarrow\\
\dot{Z}(t)=-[\dot{\theta}(t)\cdot R(\theta -\frac{\pi}{2})^{\prime }\cdot R(\theta )-\rho \cdot \Delta]\cdot Z(t)
\end{gather*}

Using the multiplicative property on rotation matrices:

\begin{equation*}
R(\theta -\frac{\pi }{2})^{\prime }=R(\frac{\pi }{2})\cdot
R(-\theta).
\end{equation*}

this leads:

\begin{equation*}
\dot{Z}(t)=-[\dot{\theta}(t)\cdot R(\frac{\pi}{2})\cdot R(-\theta )\cdot R(\theta )-\rho \cdot \Delta] \cdot Z(t)
\end{equation*}

Also:

\begin{equation*}
R(-\theta )\cdot R(\theta )=\left[
\begin{array}{cc}
1 & 0 \\
0 & 1
\end{array}
\right]
\end{equation*}

then:

\begin{equation*}
\dot{Z}(t)=-[\dot{\theta}(t)\cdot R(\frac{\pi}{2})-\rho \cdot \Delta] \cdot Z(t).
\end{equation*}

This can be recast follows:

\begin{equation*}
\dot{Z}(t)=-[\dot{\theta}(t)\cdot \left[
\begin{array}{cc}
0 & -1 \\
1 & 0
\end{array}
\right] -\rho \cdot \left[
\begin{array}{cc}
1 & 0 \\
0 & 0
\end{array}
\right]] \cdot Z(t)
\end{equation*}

In summary:

\begin{equation}\label{Sistema de ODEs originales}
\begin{cases}
\dot{X}(t)=\rho \cdot R(\theta )\cdot \Delta \cdot \underset{Z(t)}{\underbrace{%
R(\theta )^{\prime }\cdot X(t)}}\\
\dot{Z}(t)=-\left[
\begin{array}{cc}
-\rho & -\dot{\theta}(t) \\
\dot{\theta}(t) & 0
\end{array}%
\right] \cdot Z(t)\\
Z(t)=\left[
\begin{array}{cc}
z_{1}(t) & z_{2}(t)
\end{array}%
\right] ^{\prime}
\end{cases}
\end{equation}

This first order set of Ordinary Differential Equations (ODE) can be decoupled into two second order ODE's as follows:

\begin{equation*}
\begin{cases}
\ddot{z}_{1}(t)=(\rho+\frac{\ddot{\theta}(t)}{\dot{\theta}(t)}) \cdot \dot{z}_{1}(t)+(-\rho \cdot \frac{\ddot{\theta}(t)}{\dot{\theta}(t)}-\dot{\theta}(t)^{2})\cdot z_{1}(t)\\
z_{2}(t)=\frac{\dot{z_{1}(t)}-\rho \cdot z_{1}(t)}{\dot{\theta(t)}}
\end{cases}
\end{equation*}

Moreover, with the Kinematic Model for $\theta(t)$ in equation (\ref{Separacion de las dinamicas}), this set of ODE's yields:
\begin{equation}\label{ODEs de segundo orden}
\begin{cases}
\ddot{z}_{1}(t)=2 \cdot \rho \cdot \dot{z}_{1}(t)-\rho^{2} \cdot (1+\dot{\theta(t)}^{2})\cdot z_{1}(t)\\
z_{2}(t)=\frac{\dot{z_{1}(t)}-\rho \cdot z_{1}(t)}{\dot{\theta(t)}}\\
\dot{\theta(t)}=\rho \cdot \theta(t),\quad \rho<0
\end{cases}
\end{equation}

In the Appendix it is shown that (\ref{ODEs de segundo orden}) can be converted to a second order Bessel's equation which posses general solutions $J_{p}(x)$ and $Y_{p}(x)$ called Bessel's functions of first and second kind respectively -see \cite{Andrews85}, pp.228-230.

From now on and for the sake of clarity, it will used $\rho=- 1$ - otherwise the time change: $\tau=\rho \cdot t$ transform the problem to this form. In this way and using the solutions in the Appendix, the closed-form expressions for $z_{1}(t)$ and $z_{2}(t)$ using the Kinematic Model for $\theta$ in equation (\ref{Separacion de las dinamicas}), arise:

\begin{equation}\label{Soluciones en forma cerrada}
\begin{cases}
z_{1}(t)=\theta \cdot [C_{1}\cdot J_{0}(\theta)+C_{2}\cdot Y_{0}(\theta)]\\
\theta=\theta _{0}\cdot e^{- t}\\
z_{2}(t)=\frac{\dot{z}_{1}(t)+z_{1}(t)}{\dot{\theta}(t)}
\end{cases}
\end{equation}

where $C_{1}$, $C_{2}$ are arbitrary constants depending on the initial conditions.

\subsection{Asymptotic Analysis and Further Conclusions}

The analysis of the behavior of the closed-form solutions obtained when $t$ goes
to infinity will clarify some issues regarding asymptotic stability for the different immersed submanifolds forming the orbits.

In this way and using the formulas for small arguments presented in \cite{Andrews85}, pp. 248-249, the analysis for $z_{1}(t)$ from equation (\ref{Soluciones en forma cerrada}):

\begin{equation}
\lim_{t\rightarrow \infty }z_{1}(t)=\lim_{\theta \rightarrow
0}\theta\cdot [C_{1}\cdot J_{0}(\theta )+C_{2}\cdot Y_{0}(\theta )]=0
\end{equation}

providing that $J_{0}(\theta )$ is bounded $\forall \theta \in \Re $ and $\lim_{\theta \rightarrow 0}\theta \cdot Y_{0}(\theta )=0$.

The first conclusion is that $z_{1}(t)$ is asymptotically stable, however as it will be shown, this is not the case for $z_{2}(t)$. In this way, and using the property for Bessel's functions (see \cite{Andrews85}, pp. 199 and pp. 226):

\begin{equation*}
\begin{cases}
\frac{d}{dx}J_{n}(x)=\frac{n}{x}\cdot J_{n}(x)-J_{(n+1)}(x) \\
\frac{d}{dx}Y_{n}(x)=\frac{1}{2}\cdot (Y_{(n-1)}(x)-Y_{(n+1)}(x))
\end{cases}
\end{equation*}

In the case of equation (\ref{Soluciones en forma cerrada}), $\alpha=1$ and $n=0$, $z_{2}(t)$ yields:

\begin{equation}\label{z2 n=0}
z_{2}(t)=\theta \cdot [-C_{1} \cdot J_{(1)}(\theta)+\frac{C_{2}}{2}\cdot (Y_{(-1)}(\theta)-Y_{(1)}(\theta))].
\end{equation}

The formulas for small arguments -see \cite{Andrews85}, pp. 248-249, lead:

\begin{equation*}
\begin{cases}
\lim_{\theta \rightarrow 0} J_{n}(\theta)=(\frac{\theta}{2})^{n} \cdot  \frac{1}{\Gamma (n+1)}\\
\lim_{\theta \rightarrow 0}Y_{n}(\theta)=-(\frac{2}{\theta})^{n} \cdot\frac{\Gamma(n)}{\pi}
\end{cases}
\end{equation*}

where $\Gamma \left( .\right) $ is the Gamma function. Then if $\theta \rightarrow 0$, in equation (\ref{z2 n=0}):

\begin{equation*}
\lim_{\theta \rightarrow 0} z_{2}(t)=\lim_{\theta \rightarrow
0}\quad [-C_{1}\cdot \frac{\theta^{2}}{2 \cdot \Gamma(2)}-C_{2}\cdot \frac{\theta^{2} \cdot \Gamma(-1)}{4 \cdot \pi}+ C_{2}\cdot \frac{\Gamma(1)}{\pi}
\end{equation*}

Since $\Gamma(1)=\Gamma(2)=1$ and $\Gamma(-1)$ is finite, then:

\begin{equation*}
\lim_{\theta \rightarrow 0} z_{2}(t)=C_{2}\cdot \frac{1}{\pi}
\end{equation*}

To investigate the possibility for $C_{2}=0$ endowing $z_{2}$ with asymptotic stability, Definition \ref{Cambio de coordenadas} is used, taking into account that $C_{1}$ and $C_{2}$ are arbitrary constants depending on the initial conditions as pointed out in equation (\ref{Soluciones en forma cerrada}):

\begin{equation*}
\begin{array}{l}
X(t=0)=R(\theta(0))\cdot\\
\cdot
\begin{bmatrix}
\theta(0)\cdot J_{0}(\theta(0)) & \theta(0)\cdot Y_{0}(\theta(0))\\
-\theta(0)\cdot J_{0}(\theta(0)) & \theta(0)\cdot \frac{1}{2} \cdot (Y_{-1}(\theta(0))-Y_{1}(\theta(0)))
\end{bmatrix}\cdot
\begin{bmatrix}
C_{1}\\
C_{2}
\end{bmatrix}
\end{array}
\end{equation*}

where the closed-form Kinematic Model for $z_{2}$ form equation (\ref{z2 n=0}) was used. In this way, if $C_{2}=0$ is settled:

\begin{equation*}
\begin{array}{l}
X(t=0)=\theta(0) \cdot C_{1} \cdot
J_{0} \cdot R(\theta(0))\cdot
\begin{bmatrix}
1\\
-1
\end{bmatrix}
\end{array},
\end{equation*}

which is not allowing any initial condition as required to have $\Re^{2}$ a the set of attractors -see \cite{Miller82} for a complete discussion on attractors. This fact agree with the Brockett's condition regarding the impossibility to stabilize the unicycle by means of continuous controllers -see \cite{Brockett83}.

A last important conclusion arises if the set of attractors is considered:

\begin{equation*}
\begin{array}{l}
X(t\rightarrow \infty)=\underbrace{R(0)}_{I} \cdot
\begin{bmatrix}
0\\
\frac{C_{2}}{\pi}
\end{bmatrix}
\end{array}
\end{equation*}

That means that the set of attractors agree with the whole $y_{c}$ axis. On the other hand, if $\rho =1$, then applying the formulas for $\theta \rightarrow \infty$ , it can be seen that $z_{1}$ and $z_{2}$ tend asymptotically to zero leading instability in $\theta(t)$.

Next section utilizes this conclusion, extending this controller adding a dissipative term, endowing the system with an asymptotic behavior.

\section{Comparison with other Methodologies}

Once the methodologies in this paper have been presented, it is important to compare them with some of the available in the literature. In fact, the survey in \cite{Lizarraga02} will be followed in order to classify them.

\subsection{Brockett's Theorem}

As pointed out by Brockett in his theorem in \cite{Brockett83}, systems like the Unicycle can not be asymptotic stabilized using smooth controllers, however even when it is simple to arrive to this conclusion using the Brockett's theorem it is not simple to expose in a clear way the impact of this condition for every particular controller designed for the Unicycle.

In this sense the methodologies in this paper shows clearly the impossibility in using smooth controllers in order to produce asymptotic stabilization for the unicycle since closed-form solutions are available and moreover the asymptotic analysis of previous section determine this condition in an alternative way.

\subsection{Comparison with Smooth time-varying techniques}

While Method 1 and 2 in \cite{Lizarraga02}, pp. 5-13 need the evaluation of Lyapunov functions to construct stabilizing controllers for Driftless systems, the techniques in this paper is systematic and is no based in the Lyapunov theory. However when compared with the methodology of Pomet-Coron -see \cite{Lizarraga02}, pp. 14- the method of the present paper serves as a start point for the three step algorithm in \cite{Lizarraga02}, pp. 16, specially step 1.

\subsection{Comparison with Homogeneous time-varying Feedback}

Even when the theory under this subject is well developed and understood -see \cite{Lizarraga02}, pp. 20-28- the theorems in the present paper confirm some results developed using \textit{Homogeneous time-varying Feedback}, in particular compare Theorem 10 in \cite{Lizarraga02}, pp. 20 and Condition \ref{Condition for asymptotic stability} in the present paper.

\subsection{Comparison with the Transverse functions technique}

This technique is similar to the heart of the methodologies presented in this paper, in fact the stability Theorem \ref{Controlador} can be thought as a Transverse function or in other words, can be depicted as the control system running perpendicular to some Submanifolds just like the \textit{Transverse function technique}.

In this sense the Transverse function technique is more flexible and better developed but the present paper addresses a closed-form expression allowing a further analysis of the stability properties for the Unicycle.

Clearly the lack of flexibility omnipresent in the technique of the present paper is a disadvantage.

\section{The complete solution}

A slightly different model than the one consider so far can be also integrated to obtain closed-form solutions, let's consider the case presented in \cite{DelaCruz06} and write the control vector as: $u_{1}, a \cdot u_{2}$, where the parameter $a$ is the distance from $(x,y)$ to the center wheel base-line.

The, utilizing the same procedure as before:

\begin{eqnarray*}
X(t)=R(\theta) \cdot [I-R(\frac{\pi}{2})-R(\theta(t))\cdot R(-\theta(0))\cdot R(\frac{\pi}{2})] \cdot R(\theta(0))' \cdot X(0)+\\
R(\theta(t)) \cdot \int_{0}^{t} [I-R(\frac{\pi}{2})-R(\theta(t))\cdot R(-\theta(\sigma))\cdot R(\frac{\pi}{2})] \cdot u(\sigma) \cdot d\sigma
\end{eqnarray*}

where $X(t)=[x(t),y(t)]',\quad u(t)=[u_{1}, a \cdot u_{2}]'$, $'$. In this way, it is possible to analyze the effect of the parameter $a \rightarrow 0$ if $lim_{t \rightarrow \infty} \theta(t)=0$ :

\begin{eqnarray*}
lim_{t \rightarrow \infty} X(t)=\left[
\begin{bmatrix}
1 & -1\\
1 &1
\end{bmatrix}-R(-\theta(0)) \cdot R(\frac{\pi}{2})\right] \cdot X(0)+\\
\int_{0}^{\infty} \left[
\begin{bmatrix}
1 & -1\\
1 &1
\end{bmatrix}-R(-\theta(\sigma)) \cdot R(\frac{\pi}{2})\right] \cdot u(\sigma) \cdot d\sigma
\end{eqnarray*}

Clearly, when $a=0$ one control input is lost in rtdert to steer the vectro$X$ to zero. This explains why the Brockett's condition is working. It is also worth to notice that this case is of more general application (since $a=0$ include the previous case non-holonomic robots) and any control strategy like model predictive control, for instance , can be applied.

\section{Conclusions and Future Work}

A methodology to stabilize Driftless Nonlinear Control systems was presented, by applying this controller to the well known Unicycle-like vehicles, closed-form solutions for the trajectories were obtained via the solution of the second order Bessel's ODE.

The analysis of the closed-form solutions showed that the coordinate $x_{c}$ is asymptotically stable while $y_{c}$ does not. However, if the simple stabilizing law for $\theta$: $\dot{\theta}(t)=\rho\cdot \theta(t),\quad \rho<0$, which arises with this technique is considered in the opposite side, i.e. $\rho>0$, then both coordinates $x_{c},\quad y_{c}$ turn to be asymptotically stable.

This conclusion can be exploited in tow main ways:

\begin{itemize}
\item{Steer the vehicle using $\rho>0$ until the coordinates $x_{c},\quad y_{c}$ are closed enough to the origin and then to switch to $\rho<0$ stabilizing the vehicle.}
\item{Adding a dissipative term making the controller time-varying such that the closed-loop system is asymptotically stable}
\end{itemize}

As future work, several directions can be depicted, in fact, since the same constant $\rho$ is used in both components of the input control vector for the Unicycle, it would be of interest to relax this study utilizing two different constants. On the other hand, this closed-form solutions could be exploited for formations of robots where the incorporation of constraints should be considered.

\section{\label{Apendice}Appendix}

This appendix contains the proof of all the results along the paper.

\begin{proof}{\emph{(Lemma \ref{Barbalat's lemma extension})}}\\
The first step is to write $\lim_{t\rightarrow \infty}\int\limits_{0}^{t}\varphi (\tau)\cdot d\tau $ as an infinite sum:

\begin{equation*}
\lim_{t\rightarrow \infty} \int\limits_{0}^{t}\varphi (\tau)\cdot
d\tau =\sum_{i=1}^{\infty} \int\limits_{(i-1)\cdot T}^{i\cdot T}\varphi (\tau)\cdot d\tau ,
\end{equation*}

where $T \neq 0$. As it is well known for infinite series -see \cite{Knopp47}, a necessary condition for convergence is the general term tending to zero:

\begin{equation*}
\lim_{i\rightarrow \infty }\int\limits_{(i-1)\cdot T}^{i\cdot
T}\varphi (\tau )\cdot d\tau =0 ,
\end{equation*}

Since $T$ is finite and $\varphi $ is positive, then the only possibility is $\varphi(t)\rightarrow 0$ as $t\rightarrow \infty $ (or $i\rightarrow \infty $). This completes the proof.
\end{proof}


\begin{proof}{\emph{(Lemma \ref{Nuestra definicion de estabilidad})}}\\

Let first consider that the quantity $E_{c}^{*}(\infty)=\int_{0}^{\infty }\langle \dot{q}(\xi ),\dot{q}(\xi )\rangle \cdot d\xi$ is bounded and prove that the system is stable. To do this, it is required to prove that the function $E_{c}^{*}(t)=\int_{0}^{t}\left\langle \dot{q} (\xi ),\dot{q}(\xi)\right\rangle \cdot d\xi$ is strictly increasing and this can be done with the following useful result:

\vspace{2mm}

\textbf{\textit{Given a function $f:\Re\rightarrow \Re ^{+}$, then, $\varphi(t)=\int_{o}^{t}f(\tau)\cdot d\tau $ is strictly increasing.}}

\vspace{2mm}

The proof of this assertion is carried out by showing that $\int_{o}^{t_{1}}f(\tau )\cdot d\tau < \int_{o}^{t_{2}}f(\tau )\cdot d\tau $ if $t_{1}\leq t_{2}$. In fact, taking into account $\int_{t_{1}}^{t_{2}}f(\tau )\cdot d\tau >0$, then $\int_{t_{1}}^{t_{2}}f(\tau )\cdot d\tau +\int_{0}^{t_{1}}f(\tau )\cdot d\tau >\int_{0}^{t_{1}}f(\tau )\cdot d\tau $.

Then, the following Lyapunov function is chosen:

\begin{equation*}
v(q,t)=\frac{1}{E_{c}^{\ast }(t)}-\frac{1}{E_{c}^{\ast }(\infty
)}\text{.}
\end{equation*}

Realize that this function classify as a Lyapunov candidate function -see \cite{Sastry99}, pp. 188-189. In this way, we have:

\begin{equation*}
\begin{cases}
v(q=0,t)=v(0,\infty)=0,\quad \text{Since $q=0$ is an equilibrium point.}\\
v(q,t)\quad \text{is decrescent}\\
v(q,t)\quad \geq 0\\
\dot{v}(q,t)=-\frac{1}{E_{c}^{\ast }(t)^{2}}\cdot \dot{E}_{c}^{\ast }(t).
\end{cases}
\end{equation*}

Since $\dot{E}_{c}^{\ast }(t)\geq 0$ then $-\dot{v}(q,t)$ $\geq 0$
which completes the proof for this part. The counterpart of the theorem is
obvious taking into consideration that only smooth vector fields are considered.
\end{proof}


\begin{proof}{\emph{(Theorem \ref{Controlador})}}\\
The proof is based in Lemma \ref{Nuestra definicion de estabilidad} showing that the Pseudo-Kinetic Energy, $E_{c}^{*}(t)$ is bounded.

Previous the calculation of $E_{c}^{*}(t)$, it is important to notice:

\begin{equation*}
\begin{cases}
u_{i}(q)=-\rho \cdot q' \cdot S_{i}\\
i=1,2,\ldots,k\\
S_{i}^{'}\cdot S_{j}=0, \quad i\neq j\\
S_{i}^{'}\cdot S_{i}=1\\
\dot{q(t)}=\sum^{k}_{i=1} S_{i} \cdot u_{i}
\end{cases}\Rightarrow
\dot{q(t)}=\rho \cdot (\sum^{k}_{i=1} S_{i}'\cdot S_{i}) \cdot q(t)
\end{equation*}

Then, $E_{c}^{*}(t)$ arise:

\begin{equation*}
E_{c}^{*}(t)=\rho \cdot \int_{0}^{t} q(\xi)'\cdot \underbrace{(\sum^{k}_{i=1} S_{i} \cdot u_{i})}_{\dot{q(\xi)}} \cdot d\xi=\rho \cdot (\|q(t)\|^{2}-\|q(0)\|^{2})
\end{equation*}

which is bounded provided $\rho <0$.
\end{proof}


\subsection{Transforming the second order ODE in Section \ref{Trayectorias en forma cerrada},equation (\ref{ODEs de segundo orden})}

The focus in this digression is to transform the second order ODE in equation (\ref{Trayectorias en forma cerrada}) into a Bessel's ODE. To do this, the following transformation inspired in the technique presented in \cite{Andrews85}, pp. 229-230, is considered:

\begin{equation*}
z_{1}(t)=x(t)^{\alpha }\cdot l(x(t))
\end{equation*}

Then, $z_{1}$ and $\dot{z_{1}}(t)$ lead:

\begin{equation*}
\begin{cases}
\dot{z}_{1}(t)=[\alpha \cdot x(t)^{\alpha -1}\cdot l(x(t))+x(t)^{\alpha }\cdot l^{\prime}(x(t))]\cdot \dot{x}(t)\\
\ddot{z}_{1}(t)=[\alpha \cdot x(t)^{\alpha -1}\cdot l(x(t))+x(t)^{\alpha }\cdot l^{\prime}(x(t))] \cdot \ddot{x}(t)+\\
+[\alpha \cdot (\alpha -1)\cdot x(t)^{\alpha -2}\cdot l(x(t))+2\cdot \alpha \cdot x(t)^{\alpha -1}\cdot l^{\prime}(x(t))+\\
+x(t)^{\alpha }\cdot l^{\prime \prime }(x(t))]\cdot \dot{x}(t)^{2}
\end{cases}
\end{equation*}

where $\frac{dl(x)}{dx}=l^{\prime }(x(t))$ and $\frac{d^{2}l(x)}{dx^{2}}%
=l^{\prime \prime }(x(t))$. Replacing this new variables into the second order ODE at (\ref{ODEs de segundo orden}):

\begin{equation*}
\begin{array}{l}
\ddot{x}(t)\cdot [\alpha \cdot x(t)^{\alpha -1}\cdot
l(x(t))+x(t)^{\alpha }\cdot l^{\prime }(x(t))]+\\
+\dot{x}(t)^{2}\cdot \lbrack \alpha \cdot (\alpha -1)\cdot x(t)^{\alpha -2}\cdot l(x(t))+\\
2\cdot \alpha \cdot x(t)^{\alpha-1}\cdot l^{\prime }(x(t))+x(t)^{\alpha }\cdot l^{\prime\prime }(x(t))]=\\
2 \cdot \rho \cdot \dot{x}(t)\cdot \lbrack \alpha \cdot
x(t)^{\alpha -1}\cdot l(x(t))+x(t)^{\alpha }\cdot l^{\prime}(x(t)]+\\
-\rho^{2}\cdot(1+\dot{\theta}(t)^{2})\cdot x^{\alpha }\cdot l(x(t))
\end{array}
\end{equation*}

Reordering:

\begin{equation*}
\begin{array}{l}
\dot{x}(t)^{2}\cdot x(t)^{\alpha }\cdot l^{\prime \prime
}(x(t))+l^{\prime}(x(t))\cdot \lbrack \ddot{x}(t)\cdot x(t)^{\alpha }+\\
\alpha \cdot \dot{x}(t)^{2}\cdot x(t)^{\alpha -1}-2 \cdot \rho \cdot \dot{x}(t)\cdot x(t)^{\alpha}]+\\
+l(x(t))\cdot \lbrack \ddot{x}(t)\cdot \alpha \cdot x(t)^{\alpha-1}+\\
+\dot{x}(t)^{2}\cdot \alpha \cdot (\alpha -1) \cdot x(t)^{\alpha -2}-2 \cdot \rho \cdot \dot{x}(t)\cdot \alpha \cdot x(t)^{\alpha -1}+\\
+\rho^{2}\cdot(1+\dot{\theta}(t)^{2}) \cdot x(t)^{\alpha }]=0
\end{array}
\end{equation*}

In order to transform this equation into a second order \ Bessel's
ODE to the form $x^{2}\cdot \frac{d^{2}l(x)}{dx^{2}}+x\cdot \frac{dl(x)}{dx}+\left(
x^{2}-n^{2}\right) \cdot l(x)=0$, the following equalities should hold:


\begin{gather}
\vspace{2mm}
\dot{x}(t)^{2}\cdot x(t)^{\alpha }=\frac{x(t)^{2}}{f(x(t))} \label{Eq1}\\
\ddot{x}(t)\cdot x(t)^{\alpha }+2\cdot \alpha \cdot\dot{x}(t)^{2}\cdot x(t)^{\alpha -1}+ \notag\\
\vspace{2mm}
-2\cdot \rho \cdot \dot{x}(t)\cdot x(t)^{\alpha }=\frac{x(t)}{f(x(t))} \label{Eq2}\\
\ddot{x}(t)\cdot \alpha \cdot x(t)^{\alpha -1}+\dot{x}(t)^{2}\cdot \alpha \cdot (\alpha -1) \cdot x(t)^{\alpha -2}+ \notag\\
-2 \cdot \rho \cdot \dot{x}(t)\cdot \alpha \cdot x(t)^{\alpha -1}+ \notag\\
+\rho^{2}\cdot(1+\dot{\theta}(t)^{2}) \cdot x(t)^{\alpha }=\frac{x(t)^{2}-n^{2}}{f(x(t))}\label{Eq3}
\end{gather}

where $f(x(t))$ is any smooth no-null real function $\forall t\in \Re $. In order to solve this algebraic system for $\alpha $, $n$, $f(x(t))$ and $x(t)$, equations (\ref{Eq2}) and (\ref{Eq3}) are used equaling for the term $\ddot{x}(t)\cdot x(t)^{\alpha }$:

\begin{equation*}
\left\{
\begin{array}{l}
2\cdot \alpha \cdot \frac{x(t)}{f(x(t))}-2 \cdot \rho \cdot \dot{x}(t)\cdot x(t)^{\alpha}-\frac{x(t)}{f(x(t))}=\\
=\frac{x(t)\cdot (\alpha -1)}{f(x(t))}-2\cdot \rho \cdot \dot{x}(t)\cdot x(t)^{\alpha }+ \\
+\rho^{2} \cdot (1+\dot{\theta}(t)^{2})\cdot x(t)^{\alpha }\cdot \frac{x(t)}{\alpha }-\frac{(x(t)^{2}-n^{2})}{f(x(t))}\cdot \frac{x(t)}{\alpha}
\end{array}
\right.
\end{equation*}

that is:

\begin{equation*}
\alpha \cdot x(t)=\rho^{2} \cdot (1+\dot{\theta}(t)^{2})\cdot x(t)^{\alpha }\cdot \frac{x(t)}{\alpha}\cdot f(x(t))-(x(t)^{2}-n^{2})\cdot \frac{x(t)}{\alpha}
\end{equation*}

Guessing a possible solution for $x(t)$: $x(t)=\theta (t)=\theta_{0}\cdot e^{\rho\cdot t}$ with $\rho=\pm 1$:

\begin{equation*}
\begin{cases}
\alpha =x(t)^{\alpha }\cdot \frac{1}{\alpha }\cdot
f(x(t))+n^{2}\cdot \frac{1}{\alpha}\\
f(x(t))=\frac{(\alpha -\frac{n^{2}}{\alpha}) \cdot \alpha}{x(t)^{\alpha}}
\end{cases}
\end{equation*}

Replacing that function $f(x(t))$ into equation (\ref{Eq1}):

\begin{equation}\label{alpha y n}
\begin{cases}
\dot{x}(t)^{2}\cdot x(t)^{\alpha }=\frac{x(t)^{2}}{f(x(t))} \\
x(t)=\theta (t)=\theta _{0}\cdot e^{\pm t} \\
f(x(t))=\frac{(\alpha-\frac{n^{2}}{\alpha}) \cdot \alpha}{x(t)^{\alpha}}
\end{cases} \Rightarrow \quad
\alpha^{2}-n^{2}=1
\end{equation}

Using equation (\ref{Eq2}):

\begin{equation}\label{Relacion entre k1, alpha y n}
-1+2\cdot \alpha=\frac{1}{(\alpha -\frac{n^{2}}{\alpha }) \cdot \alpha}
\end{equation}

Replacing (\ref{alpha y n}) into (\ref{Relacion entre k1, alpha y n}):

\begin{equation*}
\begin{cases}
\alpha =1\\
n=\pm \sqrt{\alpha^{2}-1}=0
\end{cases}
\end{equation*}

In summary:

\begin{equation*}
\left\{
\begin{array}{l}
\rho=\pm 1\\
x(t)=\theta (t)=\theta _{0}\cdot e^{\rho \cdot t}\\
\alpha=1\\
n=0\\
z_{1}(t)=x(t)^{\alpha }\cdot l(x(t)) \\
x^{2}\cdot \frac{d^{2}l(x)}{dx^{2}}+x\cdot \frac{dl(x)}{dx}+(x^{2}-n^{2}) \cdot l(x)=0
\end{array}
\right.
\end{equation*}


\begin{thebibliography}{99}


\bibitem{Lizarraga02}  D. A. Liz\'{a}rraga, "Control of uderactuated mechanical systems using time-varying feedback and related tecnhiques," DISC Summer School 2002, Dutch Institute of Systems and Control Zeist, Netherlands, 2002.

\bibitem{Lin96} W. Lin, "Stabilization of driftless nonlinear systems via time-varying
state-feedback", Proc. Decision and Control,  1996, Korea-Japan, pp. 1803-1808.

\bibitem{Duleba98} I. Duleba and J- Sowka, "Fourier series based method of generating
continuous controls for driftless systems", International Conference on Robotics and Automation, 1998, Leuven-Belgium, pp. 3567-3572.

\bibitem{Aicardi94} M.Aicardi and G.Casalino and A.Balestrino and A. Bicchi, "Closed loop smooth stering of unicycle-like vehicles", Proceedings of the $33^{rd}$ Conference on Decision and Control, 1994, Lake Buena Vista, FL.


\bibitem{DelaCruz06} C. De la Cruz and R. Carelli, "Dynamic modeling and centralized formation control of mobile robots", Proceedings of the $32^{nd}$ of th IEEE Industrial Electronics conference, 2006.


\bibitem{Morin03} P. Morin and C. Samson, "Practical Stabilization of driftless systems on
Lie Groups: The transverse  function approach", IEEE Transactions on Automatic Control, vol. 9, 2003, pp.1496-1508.

\bibitem{Murray93} R. Murray and S. Sastry, "Nonholonomic motion planning: Steering using
sinusoids", IEEE Transactions on Autmatic Control, vol. 5, 1993, pp. 700-716.


\bibitem{Byrnes95} C. Byrnes and C. Martin,"An Integral-Invariance Principle for Nonlinear Systems", IEEE Transactions on Autmatic Control, vol. 40, 1995, pp. 983-994.




\bibitem{Brockett83}  R. W. Brockett, Asymptotic Stability and feedback stabilization, Differential Geometric Control Theory (R. W. Brockett, R. S. Millman and H. J. Sussmann, eds.), Progress in Mathematics, vol. 27. Brikh\"{a}user, 1983, pp.182-191.


\bibitem{Sastry99} S. Sastry, Nonlinear Systems: Analysis, Stability and Control,
Springer Verlag, New York, Inc., 1999.

\bibitem{Rugh96} W. J. Rugh, Linear Systems Theory, Prentice Hall, Inc., 1996.

\bibitem{Andrews85} L. C. Andrews, Special Functions for Engineers and Applied Mathematics,
Macmillan Publishing Company, 1985.

\bibitem{Miller82} R. K.Miller and A. N. Michel, Ordinary Differential Equations,
Academics Press, Inc., 1982.

\bibitem{Eastham89} M. Eastham, The Asymptotic Solution of Linear Differential Systems,
Applications of the Levinson Theorem, London Mathematical Socity Monographs, New Series, Oxford University Press, New York., 1989.

\bibitem{Horn85} R. A. Horn and C. R. Johnson, Matrix Analysis, Cambridge University
Press, 1985.

\bibitem{Knopp47} K. Knopp, Theory and Application of Infinite Series, Blackie and Son Limited, London and Glasgow, 1947, third edition (English).

\bibitem{Brockett83} R.W. Brockett and R. S. Millman and H. S. Susmann, H. S., "Asymptotic stability and feedback stabilization", Theory. Birkhäuser, Boston-Basel-Stuttgart, 1983.

\end{thebibliography}
\end{document}